\documentclass[final]{siamltex}

\usepackage{amsfonts}

\def\Frac#1#2{\frac{\displaystyle{#1}}{\displaystyle{#2}}}

\newtheorem{remark}{Remark}
\usepackage{epsfig}
\usepackage{amssymb}
\usepackage{amsmath}

\def\protectbold#1{\protect{\boldmath{$#1$}}}
\def\eqref#1{(\ref{#1})}
\def\dsp{\displaystyle}

\numberwithin{equation}{section}

\def\binomial#1#2{
\renewcommand{\arraystretch}{1.0}
\left(
\begin{array}{c} 
\hskip-5pt#1\\
\hskip-5pt#2
\end{array}
\hskip-5pt\right)}

\def\eps{\varepsilon}

\makeatother

\def\dsp{\displaystyle}

\def\FG#1#2#3#4{
{}_2F_1\left(
\begin{array}{c}
\begin{array}{c}\hskip-10pt#1,#2\end{array}\\
\begin{array}{c}\hskip-10pt #3\end{array}
\end{array}
\hskip-8pt;\,#4
\right)}

\def\bigO{{\cal O}}
\def\calC{{{\cal C}}}

\title{Non-iterative computation of Gauss-Jacobi quadrature \thanks{
This work was supported by  {\it Ministerio de Econom\'{\i}a y Competitividad}, 
project MTM2015-67142-P (MINECO/FEDER, UE). NMT thanks CWI, Amsterdam, for scientific support.
}}

\author{Amparo Gil\thanks{
Departamento de Matem\'atica Aplicada y CC. de la Computaci\'on.
ETSI Caminos. Universidad de Cantabria. 39005-Santander, Spain.
({\tt amparo.gil@unican.es}).}
 \and Javier Segura \thanks{Departamento de Matem\'aticas, Estad\'{\i}stica y Computaci\'on.
Facultad de Ciencias. Universidad de Cantabria. 39005-Santander, Spain.
({\tt javier.segura@unican.es}).}
        \and Nico M. Temme\thanks{IAA, 1825 BD 25, Alkmaar, The Netherlands.
{Former address: Centrum Wiskunde \& Informatica (CWI), 
Science Park 123, 1098 XG Amsterdam,  The Netherlands}
({\tt nico.temme@cwi.nl}).} 
}

\begin{document}

\maketitle

\begin{abstract}
Asymptotic approximations to the zeros of Jacobi polynomials are given, with methods to obtain the coefficients in the expansions. These approximations 
can be used as standalone methods for the non-iterative computation of the nodes of 
Gauss--Jacobi quadratures of high degree ($n\ge 100$). 
We also provide asymptotic approximations for functions related to the first order derivative of Jacobi polynomials
which are used for computing  the weights of the Gauss--Jacobi quadrature.
The performance of the asymptotic approximations is
illustrated with numerical examples, and it is shown that nearly double precision relative accuracy 
is obtained both for the nodes and the weights when $n\ge 100$ and $-1< \alpha, \beta\le 5$. 
For smaller degrees the approximations are also useful as they provide $10^{-12}$ relative accuracy for 
the nodes when $n\ge 20$, and just one Newton step would be sufficient to guarantee double precision accuracy 
in that cases.
\end{abstract}

\begin{keywords} 
Jacobi polynomials; asymptotic expansions; Gaussian quadratures.
\end{keywords}

\begin{AMS}
33C45, 41A60, 65D20, 65D32
\end{AMS}

\pagestyle{myheadings}
\thispagestyle{plain}
\markboth{A. Gil, J. Segura, N.M. Temme}{Manuscript submitted to SIAM J. Sci. Comp.}

\section{Introduction}\label{sec:Intro}

Given a weight function $w(x)$ in an interval $[a,b]$, the Gaussian quadrature formula with $n$ nodes 
$Q_n (f)=\sum_{i=1}^n w_i f(x_i)$ for approximating the integral $I(f)=\int_a^b f(x) w(x) \,dx$ is the formula with 
the highest possible degree of accuracy, that is, the formula such that $I(f)=Q_n(f)$ if $f$ is any polynomial
of degree smaller than $2n$. Because of the optimality with respect to the degree of accuracy, Gaussian rules
have fast convergence as the degree increases (especially for analytic functions), and they are one of the most 
popular methods of numerical integration, appearing in countless applications. 
A particularly important set of quadrature rules are the Gauss-Jacobi rules, with weight functions
$w(x)=(1-x)^{\alpha}(1+x)^{\beta}$, $\alpha,\beta>-1$ in $[a,b]=[-1,1]$ (Gauss-Legendre quadrature is the particular case
$\alpha=\beta=0$). 

The Golub-Welsch algorithm \cite{Golub:1969:COG} is an interesting
method for computing quadrature rules of low degree. However, as the number of
nodes n increases, the complexity scales as ${\cal O}(n^2)$.
Recent papers describing fast methods of order ${\cal O}(n)$ for the computation of classical Gauss quadratures are 
\cite{Glaser:2007:AFA,Bogaert:2012:OCO,Hale:2013:FAC,Bogaert:2014:IFC,Townsend:2016:FCO,
Bremer:2017:OTN,Gil:2018:GHL}. 
A nice short review describing recent developments in the computation
of high degree Gauss-Legendre quadrature can be found in \cite{Town:2015:TRF}. 

In \cite{Hale:2013:FAC}, the method for Gauss-Jacobi quadrature 
relies on simple asymptotic approximations for the nodes which
are iteratively refined by Newton's method, with the  Jacobi polynomial (whose roots are the nodes of the 
formula) and its derivative computed by two types of available asymptotic expansions; 
the same approach was used in \cite{Townsend:2016:FCO} for Gauss-Hermite quadrature. 
A similar strategy was considered earlier for Gauss-Legendre quadrature in 
\cite{Swarztrauber:2002:OCT,Glaser:2007:AFA,Bogaert:2012:OCO}, but using different methods 
for computing the Legendre polynomials; 
methods for Gauss-Laguerre and Gauss-Hermite were also described in \cite{Glaser:2007:AFA}. 
A more recent approach, based in the integration 
of ODEs using non-oscillatory phase representations is that of \cite{Bremer:2017:OTN}, which appears to be
competitive for Gauss rules of very high degree. All these approaches have in common the use of an iterative
root-finding method (Newton's method in all cases) for polishing the initial estimations for the roots.

However, 
for the case of Gauss-Legendre quadrature, it was shown in \cite{Bogaert:2014:IFC} that the direct computation of 
the nodes and
weights by explicit asymptotic expansions was possible, without recourse to iterative methods, and that this non-iterative
approach is accurate for moderate orders ($n\ge 100$) and faster than iterative methods. In \cite{Gil:2018:GHL} it
was shown that for the other two ``classical" quadrature formulas (Gauss-Hermite and Gauss-Laguerre) a purely asymptotic 
approach was also possible. Finally, in the present paper we consider the asymptotic computation of the nodes and weights 
of Gauss-Jacobi quadrature, which completes the description of large-degree asymptotic methods for the classical Gauss quadrature rules (Hermite, Laguerre, Jacobi). 

The present paper extends the results in \cite{Bogaert:2014:IFC}, where the case of the Gauss-Legendre 
quadrature ($\alpha=\beta=0$) was considered, as well as those in \cite{Hale:2013:FAC}, which described an iterative method producing close to double precision absolute accuracy for the nodes. We obtain a fast non-iterative method based on asymptotic expansions for the Jacobi nodes and weights with nearly double precision relative accuracy for degrees $n\ge 100$ and $-1<\alpha, \beta \le 5$.

We consider asymptotic approximations of two types, one in terms of Bessel functions and a simpler one in terms of elementary functions. In \cite{Bogaert:2014:IFC}, the computation was based on Bessel-type expansions, but we show that our simpler 
elementary approach is in fact able to compute accurately most of the nodes and weights (although the Bessel expansions should
be considered for the extreme nodes). We discuss how to invert these asymptotic series in order to obtain asymptotic expansions for the nodes, and how to obtain the weights in a numerically stable way using asymptotic expansions related to the derivative of Jacobi polynomials.

We start the paper by summarizing some basic properties satisfied
by Jacobi polynomials which are relevant to the paper. In \S\ref{sec:Jacpolelem} we derive two large-$n$ asymptotic 
expansions in terms of elementary functions valid for $x\in[-1+\delta,1-\delta]$, where $\delta$ is a small positive number.  Always, $\alpha$ and $\beta$ should be of order $o(1)$ for large values of $n$. The expansion in \eqref{eq:Jacelem08}--\eqref{eq:Jacelem09} is used in \S\ref{sec:Jaczerelem} for finding asymptotic expansions of the zeros. The expansion in terms of the $J$-Bessel function is considered 
in \S\ref{sec:JacpolBess}; see \eqref{eq:JacBes04}. This result is valid for $x\in[-1+\delta,1]$, where again $\delta$ is a small positive number, and it is used in \S\ref{sec:JacBeszer} to find asymptotic expansions of the zeros.  In \S\ref{sec:numerical} we give details of the numerical performance of these expansions. Finally, in \S\ref{paltonto} we summarize other applications of the methods described 
in this paper, in particular Gauss-Lobatto quadrature and barycentric interpolation.

\subsection{Jacobi polynomials: first properties}

The Jacobi polynomials $P^{(\alpha,\beta)}_{n}(x)$, with $\alpha, \beta > -1$, constitute an orthogonal set with weight function $w(x)=(1-x)^{\alpha}(1+x)^{\beta}$ on the interval $[-1,1]$, and they satisfy the symmetry relation
\begin{equation}\label{eq:Intro01}
P_n^{(\alpha,\beta)}(-x)=(-1)^nP_n^{(\beta,\alpha)}(x).
\end{equation}
An explicit representation can be given in the form of a Gauss hypergeometric function
\begin{equation}\label{eq:Intro02}
\begin{array}{@{}r@{\;}c@{\;}l@{}}
P^{(\alpha,\beta)}_{n}(x)&=&\dsp{\sum_{\ell=0}^{n}\frac{{\left(n+\alpha+%
\beta+1\right)_{\ell}}{\left(\alpha+\ell+1\right)_{n-\ell}}}{\ell!\;(n-\ell)!}%
\left(\frac{x-1}{2}\right)^{\ell}}\\[8pt]
&=&\dsp{\frac{{\left(\alpha+1\right)_{n}}}{n!}\FG{-n}{n+\alpha+\beta+1}{\alpha+1}{\tfrac12(1-x)}},
\end{array}
\end{equation}
and the many transformation formulas of the Gauss function give many other representations.

The Rodrigues formula is
\begin{equation}\label{eq:Intro03}
P_n^{(\alpha,\beta)}(x)=\frac{(-1)^n}{2^n n!\,w(x)}\frac{d^n}{dx^n}\left(w(x)(1-x^2)^n\right),
\end{equation}
and this gives the integral representation
\begin{equation}\label{eq:Intro04}
P_n^{(\alpha,\beta)}(x)=\frac{(-1)^n}{2^n\,w(x)}\frac{1}{2\pi i}\int_\calC \frac{w(z)(1-z^2)^n}{(z-x)^{n+1}}\,dz, \quad x\in(-1,1),
\end{equation}
where the contour $\calC$ is a circle around the point $z=x$ with a radius small enough to have the points $\pm1$ outside the circle. 

The recurrence relation is
\begin{equation}\label{eq:Intro05}
P_{n+1}^{(\alpha,\beta)}(x)=(A_n x+B_n)P_n^{(\alpha,\beta)}(x)-C_nP_{n-1}^{(\alpha,\beta)}(x),
\end{equation}
where
\begin{equation}\label{eq:Intro06}
\begin{array}{@{}r@{\;}c@{\;}l@{}}
A_{n}&=&\dfrac{(2n+\alpha+\beta+1)(2n+\alpha+\beta+2)}{2(n+1)(n+\alpha+\beta+1)},\\[8pt]
B_{n}&=&\dfrac{(\alpha^{2}-\beta^{2})(2n+\alpha+\beta+1)}{2(n+1)(n+\alpha+\beta+%
1)(2n+\alpha+\beta)},\\[8pt]
C_{n}&=&\dfrac{(n+\alpha)(n+\beta)(2n+\alpha+\beta+2)}{(n+1)(n+\alpha+\beta+1)(2%
n+\alpha+\beta)},
\end{array}
\end{equation}
with initial values $P_{0}^{(\alpha,\beta)}(x)=1$, $P_{1}^{(\alpha,\beta)}(x)=\frac12(\alpha-\beta)+\frac12(\alpha+\beta+2)x$.

\section{Expansions in terms of elementary functions}\label{sec:Jacpolelem}
We give details for an asymptotic expansion of the Jacobi polynomial $P_n^{(\alpha,\beta)}(x)$ that is valid inside the interval $[-1+\delta,1-\delta]$, with $\delta$ a fixed small positive number. This expansion is best suited 
for the zeros of  $P_n^{(\alpha,\beta)}(x)$ around the origin, but in fact it will provide double precision
 accuracy for most of the zeros, as we will discuss. Because of the relation  $P_n^{(\alpha,\beta)}(x)= (-1)^nP_n^{(\beta,\alpha)}(-x)$ we can concentrate on the positive zeros. 

\subsection{An expansion derived by E. Hahn}\label{sec:Hahn}
In \cite[\S18.15(i)]{Koornwinder:2010:OPS} an expansion is given derived in \cite{Hahn:1980:AJP}, which has the nice property that the coefficients are known in explicit form. 

The expansion is described in terms of several formulas. We have for large values of $n$
\begin{equation}\label{eq:Jacelem01}
\begin{array}{@{}r@{\;}c@{\;}l@{}}
P^{(\alpha,\beta)}_{n}\left(\cos\theta\right)\ %
&=&\dsp{\frac{2^{2n+\alpha+\beta+1}\mathrm{B}\left(n+\alpha+1,n+\beta+1\right)}
{\pi\,\sin^{\alpha+\frac{1}{2}}\frac{1}{2}\theta \,
\cos^{\beta+\frac{1}{2}}\frac{1}{2}\theta}\ \times}\\[8pt]
&&
\dsp{\left(\sum_{m=0}^{M-1}\frac{f_{m}(\theta)}{2^{m}{\left(2n+\alpha+\beta+2\right%
)_{m}}}+\bigO\left(n^{-M}\right)\right),}
\end{array}
\end{equation}
where
\begin{equation}\label{eq:Jacelem02}
\begin{array}{@{}r@{\;}c@{\;}l@{}}
f_{m}(\theta)&=&\dsp{\sum_{\ell=0}^{m}\frac{C_{m,\ell}(\alpha,\beta)}{\ell!(m-\ell)!}%
\frac{\cos\theta_{n,m,\ell}}{\left(\sin\frac{1}{2}\theta\right)^{\ell}\left(%
\cos\frac{1}{2}\theta\right)^{m-\ell}},}\\[12pt]
C_{m,\ell}(\alpha,\beta)&=&{\left(\frac{1}{2}+\alpha\right)_{\ell}}{\left(%
\frac{1}{2}-\alpha\right)_{\ell}}{\left(\frac{1}{2}+\beta\right)_{m-\ell}}{%
\left(\frac{1}{2}-\beta\right)_{m-\ell}},\\[12pt]
\theta_{n,m,\ell}&=&\frac{1}{2}(2n+\alpha+\beta+m+1)\theta-\frac{1}{2}(\alpha+%
\ell+\frac{1}{2})\pi.
\end{array}
\end{equation}

This expansion has been used in a paper by Gatteschi and Pittaluga
\cite{Gatteschi:1985:AAE} to derive an approximation of the mid range zeros, and this
approximation in terms of elementary functions 
is used in \cite{Hale:2013:FAC} as first estimates for these zeros. Additionally, in \cite{Hale:2013:FAC}
the expansion (\ref{eq:Jacelem01}) (and that for its derivative) is used for the iterative refinement
of the first estimations of these zeros.
 
 For deriving more accurate approximations of the mid-zeros we prefer a different expansion, 
which will be given in the next section. Differently to 
\cite{Hale:2013:FAC}, iterative methods (Newton) will not be needed in order to improve the accuracy; 
after a node has been computed by means of the asymptotic expansion, we will only need to evaluate 
one additional asymptotic expansion in order to compute the corresponding weight. We therefore estimate
that our computations will be around factor $2N$ faster that the iterative algorithm \cite{Hale:2013:FAC}, 
with $N$ the average number of iterations used in 
the algorithm of \cite{Hale:2013:FAC} ($2N$ and not $N$ because the expansions both for the polynomial
and the derivative are used in the iterative algorithm).

\subsection{A compound Poincar\'e-type expansion}\label{sec:Jacpelem}
In this section we give an expansion which has the canonical form
\begin{equation}\label{eq:Jacelem03}
P_n^{(\alpha,\beta)}(\cos\theta)=\frac{1}{\sqrt{\pi \kappa}}\frac{\cos \chi \,P(x)-\sin\chi \,Q(x)}{\sin^{\alpha+\frac12}\frac12\theta\,\cos^{\beta+\frac12}\frac12\theta},
\end{equation}
where
\begin{equation}\label{eq:Jacelem04}
x=\cos\theta,\quad \chi=\kappa\theta-\left(\tfrac12\alpha+\tfrac14\right)\pi,\quad \kappa = n+\tfrac12(\alpha+\beta+1),
\end{equation}
with expansions
\begin{equation}\label{eq:Jacelem05}
P(x)\sim\sum_{m=0}^\infty \frac{p_m(x)}{\kappa^m},\quad Q(x)\sim\sum_{m=0}^\infty \frac{q_m(x)}{\kappa^m},\quad \kappa\to\infty,
\end{equation}
for $x\in [-1+\delta,1-\delta]$, and $\alpha$ and $\beta$ bounded. A few steps to obtain these expansions will be given in \S\ref{sec:Jacsadelem}.

The first coefficients are
\begin{equation}\label{eq:Jacelem06}
\begin{array}{@{}r@{\;}c@{\;}l@{}}
p_0(x)&=&1, \quad  q_0(x)=0,\\[8pt]
p_1(x)&=&-\frac12\alpha\beta,\quad  \dsp{q_1(x)=\frac{2\alpha^2-2\beta^2+(2\alpha^2+2\beta^2-1)x}{8\sin\theta}},\\[8pt]
p_2(x)&=&-\bigl(4\alpha^4+4\beta^4-16\beta^2-16\alpha^2-24\alpha^2\beta^2+8\ +\\[8pt]
&&\quad
4(2\alpha^2+2\beta^2-5)(\alpha^2-\beta^2)x\ +\\[8pt]
&&\quad
(4\alpha^4+4\beta^4+24\alpha^2\beta^2-4\alpha^2-4\beta^2+1)x^2\bigr)/(128\sin^2\theta),\\[8pt]
q_2(x)&=&-\frac12\alpha\beta \,q_1(x).
\end{array}
\end{equation}

The expansions in \eqref{eq:Jacelem05} have negative powers of $\kappa$. We can modify the expansions by introducing the front factor\footnote{This function is also present in the expansions in terms of Bessel functions, 
see \S\ref{sec:JacpolBess}.}
\begin{equation}\label{eq:Jacelem07}
G_\kappa(\alpha,\beta)=\frac{\Gamma(n+\alpha+1)}{n!\,\kappa^{\alpha}}=\frac{\Gamma\left(\kappa+\frac12(\alpha-\beta+1)\right)}{\Gamma\left(\kappa-\frac12(\alpha+\beta-1)\right)\,\kappa^{\alpha}},
\end{equation}
and by multiplying the expansions in  \eqref{eq:Jacelem05}  by the asymptotic expansion of $1/G_\kappa(\alpha,\beta)$.
The result is the representation
\begin{equation}\label{eq:Jacelem08}
P_n^{(\alpha,\beta)}(\cos\theta)=\frac{G_\kappa(\alpha,\beta)}{\sqrt{\pi \kappa}}\frac{\cos \chi \,U(x)-\sin\chi \,V(x)}
{\sin^{\alpha+\frac12}\frac12\theta\,\cos^{\beta+\frac12}\frac12\theta},
\end{equation}
with expansions
\begin{equation}\label{eq:Jacelem09}
U(x)\sim\sum_{m=0}^\infty \frac{u_{2m}(x)}{\kappa^{2m}},\quad V(x)\sim\sum_{m=0}^\infty \frac{v_{2m+1}(x)}{\kappa^{2m+1}}.
\end{equation}
The first coefficients are
\begin{equation}\label{eq:Jacelem10}
\begin{array}{@{}r@{\;}c@{\;}l@{}}
u_0(x)&=&1, \quad  v_1(x)=q_1(x),\\[8pt]
u_2(x)&=&\dsp{\frac{1}{384\sin^2\theta}}\bigl(
12(5-2\alpha^2-2\beta^2)(\alpha^2-\beta^2)x\ +\\[8pt]
&&\quad
4(-3(\alpha^2-\beta^2)^2+3(\alpha^2+\beta^2)-6+4\alpha(\alpha^2-1+3\beta^2)\ +\\[8pt]
&&\quad
(-12(\alpha^2+\beta^2)(\alpha^2+\beta^2-1)-16\alpha(\alpha^2-1+3\beta^2)-3)x^2\bigr).
\end{array}
\end{equation}

We give the first terms of the expansion of $G_\kappa(\alpha,\beta)$:
 \begin{equation}\label{eq:Jacelem11}
\begin{array}{@{}r@{\;}c@{\;}l@{}}
G_\kappa(\alpha,\beta)
&\sim& \dsp{1-\frac{\alpha\beta}{2\kappa} -\frac{\alpha(\alpha-1)(1+\alpha-3\beta^2)}{24\kappa^2}\ +}\\[8pt]
&&
\dsp{\frac{\alpha(\alpha-1)(\alpha-2)\beta(1+\alpha-\beta^2)}{48\kappa^3}+\ldots.}
\end{array}
\end{equation}
A more efficient expansion in negative powers of $\left(\kappa-\frac12\beta\right)^2$ reads
 \begin{equation}\label{eq:Jacelem12}
G_\kappa(\alpha,\beta)\sim (w/\kappa)^{\alpha}\sum_{m=0}^\infty \frac{C_m(\rho) (-\alpha)_{2m}}{w^{2m}},
\end{equation}
where
 \begin{equation}\label{eq:Jacelem13}
w=\kappa-\tfrac12\beta, \quad \rho=\tfrac12(\alpha+1),
\end{equation}
and the first $C_m(\rho)$ are
 \begin{equation}\label{eq:Jacelem14}
C_0(\rho)=1, \quad C_1(\rho)=-\tfrac1{12}\rho, \quad C_2(\rho)=\tfrac1{1440}\left(5\rho+1\right).
\end{equation}
For more details we refer to \cite[\S5.11(iii)]{Askey:2010:GFN} and \cite[\S6.5.2-6.5.3]{Temme:2015:AMI}.

\begin{remark}\label{rem:rem01}
When $\alpha=\beta=-\frac12$ the Jacobi polynomial becomes a Chebyshev polynomial. We have
\begin{equation}\label{eq:Jacelem15}
P_n^{(-\frac12,-\frac12)}(\cos\theta)=\frac{2^{-2n}(2n)!}{(n!)^2}\cos (n\theta)=\frac{\Gamma\left(n+\frac12\right)}{\sqrt{\pi}\,\Gamma(n+1)}\cos (n\theta).
\end{equation}
In this special case we have $\kappa=n$, $\chi=n\theta$. The  $q_m(x)$ and $v_m(x)$ all vanish, $p_m(x)$ and $u_m(x)$ all become independent of $x$, and we obtain 
\mbox{$P(x)=\sqrt{n}\,\Gamma\left(n+\frac12\right)/\Gamma(n+1)$}.
\end{remark}

\subsubsection{Expansions of derivatives}\label{sec:Jacderelem}

A representation of the derivative can be written in the form
\begin{equation}\label{eq:Jacelem16}
\frac{d}{d\theta}P_n^{(\alpha,\beta)}(\cos\theta)=\sqrt{\frac{\kappa}{\pi}}G_\kappa(\alpha,\beta)\frac{\sin \chi \,Y(x)+\cos\chi \,Z(x)}{\sin^{\alpha+\frac12}\frac12\theta\,\cos^{\beta+\frac12}\frac12\theta},
\end{equation}
with expansions
\begin{equation}\label{eq:Jacelem17}
Y(x)\sim\sum_{m=0}^\infty \frac{y_{2m}(x)}{\kappa^{2m}},\quad Z(x)\sim\sum_{m=0}^\infty \frac{z_{2m+1}(x)}{\kappa^{2m+1}},
\end{equation}
where $y_0(x)=1$ and for $m=0,1,2,\ldots$
\begin{equation}\label{eq:Jacelem18}
\begin{array}{@{}r@{\;}c@{\;}l@{}}
y_{2m}(x)&=&\dsp{u_{2m}(x)+\frac{d}{d\theta} v_{2m-1}(x)+\frac{\beta-\alpha-(\alpha+\beta+1)x}{2\sin\theta}v_{2m-1}(x),}\\[8pt]
z_{2m+1}(x)&=&\dsp{\frac{d}{d\theta} u_{2m}(x)-v_{2m+1}(x)+\frac{\beta-\alpha-(\alpha+\beta+1)x}{2\sin\theta}u_{2m}(x).}
\end{array}
\end{equation}

For computing the Gauss weights it is convenient to have an expansion of the derivative of the function
\begin{equation}\label{eq:Jacelem19}
W(\theta)=\cos \chi \,U(x)-\sin\chi \,V(x),
\end{equation}
which is the oscillatory part in the representation given in \eqref{eq:Jacelem08}. We have
\begin{equation}\label{eq:Jacelem20}
\frac{d}{d\theta}W(\theta)=-\kappa\bigl(\sin\chi M(x)+\cos\chi N(x)\bigr),
\end{equation}
with expansions
\begin{equation}\label{eq:Jacelem21}
M(x)\sim \sum_{\ell=0}^\infty\frac{m_{2\ell}(x)}{\kappa^{2\ell}},\quad
N(x)\sim \sum_{\ell=0}^\infty\frac{n_{2\ell+1}(x)}{\kappa^{2\ell+1}},
\end{equation}
where $m_0(x)=u_0(x)=1$,  $n_1(x)=v_1(x)$, and for $\ell=1,2,3,\ldots$
\begin{equation}\label{eq:Jacelem22}
\begin{array}{@{}r@{\;}c@{\;}l@{}}
m_{2\ell}(x)&=&\dsp{u_{2\ell}(x)-\sin\theta \frac{d}{dx}v_{2\ell-1}(x)},\\[8pt]
n_{2\ell+1}(x)&=&\dsp{v_{2\ell+1}(x)+\sin\theta  \frac{d}{dx}u_{2\ell}(x)}.  
\end{array}
\end{equation}

\subsubsection{An expansion of the zeros}\label{sec:Jaczerelem}
We denote the zeros of $P_n^{(\alpha,\beta)}(\cos\theta)$ by $x_1,x_2,\cdots,x_n$.

As a first-order approximation in terms of the $\theta$ variable we take
\begin{equation}\label{eq:Jacelem23}
\theta_0=\frac{\left(n+1-k+\frac12\alpha-\frac14\right)\pi}{\kappa},\quad k=1,2,\ldots,n.
\end{equation}
For this value of $\theta$, $\chi$ defined in \eqref{eq:Jacelem04} becomes $\chi_0=(n-k+\frac12)\pi$, which gives $\cos\chi_0=0$. This initial value $\theta_0$ is given in \cite{Gatteschi:1985:AAE},
together with the expansion $\theta=\theta_0+\theta_1/\kappa^2+\bigO\left(\kappa^{-4}\right)$, with 
\begin{equation}\label{eq:Jacelem24}
\theta_1=\left(\tfrac{1}{16}-\tfrac14\alpha^2\right)\cot\left(\tfrac12\theta_0\right)-\left(\tfrac{1}{16}-\tfrac14\beta^2\right)\tan\left(\tfrac12\theta_0\right).
\end{equation}
This is $-v_1(x)=-q_1(x)$, see \eqref{eq:Jacelem06}.

To find a few more details, we use $W(\theta)$ defined in \eqref{eq:Jacelem19}, and expand this function at $\theta_0$ by writing $\theta=\theta_0+\eps$, which gives for a zero $\theta$
\begin{equation}\label{eq:Jacelem25}
W(\theta)=W(\theta_0+\eps)=W(\theta_0)+\frac{\eps}{1!}W^\prime(\theta_0)+ \frac{\eps^2}{2!}W^{\prime\prime}(\theta_0)+\ldots = 0,
\end{equation}
where the derivatives are with respect to $\theta$. We assume for $\eps$ an expansion in the form
\begin{equation}\label{eq:Jacelem26}
\eps\sim\frac{\theta_1}{\kappa^2}+\frac{\theta_2}{\kappa^4}+\frac{\theta_3}{\kappa^6}+\ldots.
\end{equation}
Using this expansion and those of $U(x)$ and $V(x)$ in \eqref{eq:Jacelem09}, and comparing equal powers of $\kappa$, we find 
\begin{equation}\label{eq:Jacelem27}
\begin{array}{@{}r@{\;}c@{\;}l@{}}
\theta_1&=& -v_1,\\[8pt]
\theta_2&=& u_2v_1+v_1^{\prime}v_1+\frac{1}{3}v_1^3-v_3,\\[8pt]
\theta_3&=& -\frac{4}{3}v_1^{\prime}v_1^3-\frac{1}{5}v_1^5-v_5+v_3v_1^2-\frac{1}{2}v_1^{\prime\prime}v_1^2-2v_1^{\prime}u_2v_1-(v_1^{\prime})^2v_1\ +\\[8pt]
&&\quad
v_1^{\prime}v_3-u_2^2v_1-u_2^{\prime}v_1^2+u_4v_1-u_2v_1^3+u_2v_3+v_3^{\prime}v_1,
\end{array}
\end{equation}
where $u_m,v_m$ have argument $x$, and $x=\cos(\theta_0)$. The derivatives are with respect to $\theta$.
In terms of the original variables:
 \begin{equation}\label{eq:Jacelem28}
\begin{array}{@{}r@{\;}c@{\;}l@{}}
\theta_1&=& \dsp{ -\frac{1}{8\sin\theta}}\left(2\beta^2x+2\alpha^2x-x-2\beta^2+2\alpha^2\right),\\[8pt]
\theta_2&=&   \dsp{\frac{1}{384\sin^3\theta}}\bigl(-33x-36\beta^2x^2+36\alpha^2x^2+24\beta^4x^2-24\alpha^4x^2+2x^3\ +\\[8pt]
&&\quad
84\beta^2x-60\alpha^4x-60\beta^4x+84\alpha^2x+4\beta^4x^3+4\alpha^4x^3-8\beta^2x^3\ +\\[8pt]
&&\quad
40\alpha^2-8\alpha^2x^3-40\beta^2+32\beta^4-32\alpha^4+24\alpha^2\beta^2x^3-24\alpha^2\beta^2x\bigr).
\end{array}
\end{equation}

We summarize the steps for finding a zero $x_k$, $1\le k \le n$, given $n$, $\alpha$, and $\beta$.
\begin{enumerate}
\item
Compute $\theta_0$ defined in \eqref{eq:Jacelem23}.
\item
With this $\theta_0$, compute  the coefficients given in \eqref{eq:Jacelem27} or  \eqref{eq:Jacelem28}, with $x=\cos\theta_0$. The coefficients $u_{2m}(x)$ and $v_{2m+1}(x)$ are those in the expansions in \eqref{eq:Jacelem09}.
\item
Compute $\eps$ from \eqref{eq:Jacelem26} and next $\theta=\theta_0+\eps$.
\item 
Compute $x_k\sim\cos\theta$. 
\end{enumerate}

\subsubsection{Details on computing \protectbold{\cos\chi} }\label{sec:Jaccos}
When we evaluate the functions $W(\theta)$ or $W^\prime(\theta)$ (see \eqref{eq:Jacelem19} and \eqref{eq:Jacelem20}), with $\theta$ in the form $\theta=\theta_0+\eps$, see \eqref{eq:Jacelem23}, we can write $\chi$ defined in \eqref{eq:Jacelem04} as
\begin{equation}\label{eq:Jacelem29}
\chi=\kappa\eps+\left(n-k+\tfrac12\right)\pi,
\end{equation}
and, hence,
\begin{equation}\label{eq:Jacelem30}
\cos\chi=(-1)^{n-k+1}\sin(\kappa\eps),\quad \sin\chi=(-1)^{n-k}\cos(\kappa\eps).
\end{equation}
In the original forms, especially for the zeros in the middle of the interval,  the argument $\chi$  may be of order $\kappa$, in the new forms the arguments are of order $1/\kappa$.  As a consequence,  the evaluation of these functions can be done with better accuracy when we use \eqref{eq:Jacelem30}.

For example, when we take $n= 100$,  $\alpha= \frac13$, $\beta= \frac15$ and 
the expansion in \eqref{eq:Jacelem26} with the given three terms in \eqref{eq:Jacelem27}, we obtain for the middle zero $x_{50}$, using Maple with Digits=16, 
\begin{equation}\label{eq:Jacelem31}
\begin{array}{@{}r@{\;}c@{\;}l@{}}
\cos\chi&=&0.0001908363241002135, \\
 (-1)^{n-k+1}\sin(\kappa\eps)&=&0.0001908363242842724,
\end{array}
\end{equation}
with a relative error $9.64\times10^{-10}$. Corresponding   values of $W(\theta)$ are, when using both forms,
\begin{equation}\label{eq:Jacelem32}
 -0.1837567\times 10^{-12} \quad {\rm and} \quad 0.3013\times10^{-15}.
\end{equation}
This feature is more relevant for the middle zeros than for the other ones. However, the middle zeros are more interesting than those near the endpoints  $\pm1$, because the expansions are not valid there.

A similar problem may occur when evaluating $x=\cos\theta$ when we have found the approximation of a zero $\theta$ near $\frac12\pi$, that is, near $x=0$. Assume that $n$ is even, and write $n=2m$. Then $m=\frac12\kappa-\frac14(\alpha+\beta+1)$ and for a zero $x_k$ near the origin, we write $k=m+j$. Then, $\theta_0=\frac12\pi+\tau$, where 
$\tau=(\alpha-\beta+2-4j)/(4\kappa)$, and $x=\cos(\theta_0+\eps)=-\sin(\eps+\tau)$. Similar for odd $n=2m+1$ and the zero $x_k$ with $k=m+j$. Then we have $\tau=(\alpha-\beta+4-4j)/(4\kappa)$. When $\alpha=\beta$ and $j=1$ (the zero at the origin in this case), we have $\tau=0$.
    
\subsubsection{The saddle point method for obtaining the expansion}\label{sec:Jacsadelem}
To obtain the representation in \eqref{eq:Jacelem03} and the expansions in \eqref{eq:Jacelem05} we use the integral representation given in \eqref{eq:Intro04}
where $w(z)=(1-z)^\alpha(1+z)^\beta$. We write the integral in the form
\begin{equation}\label{eq:Jacelem33}
P_n^{(\alpha,\beta)}(x)=\frac{(-1)^n}{2^nn!\,w(x)}\frac{1}{2\pi i}\int_\calC e^{-\kappa\phi(z)}\frac{w(z)(z-x)^{\gamma-1}}{(1-z^2)^{\gamma}}\,dz, 
\end{equation}
where $\kappa$ is given in \eqref{eq:Jacelem04} and 
\begin{equation}\label{eq:Jacelem34}
\gamma=\tfrac12(\alpha+\beta+1), \quad \phi(z)=\ln(z-x)-\ln(1-z^2).
\end{equation}
The saddle points $z_\pm$ follow from the zeros of
\begin{equation}\label{eq:Jacelem35}
\phi^\prime(z)=\frac{z^2-2xz+1}{(z-x)(1-z^2)},\quad \Longrightarrow\quad z_\pm=e^{\pm i \theta },\quad x=\cos\theta.
\end{equation}
The saddle point contour runs through $z_-$ from $z=-1$ to $z=1$ (below the real $z$-axis), and then through $z_+$ from $z=1$ to $z=-1$ (above the real axis).

The contribution $P^+$ from $z_+$ follows from the substitution
\begin{equation}\label{eq:Jacelem36}
\phi(z)-\phi(z_+)=\tfrac12s^2,
\end{equation}
with corresponding points $z=\pm1$ and $s=\pm\infty$. This gives 
\begin{equation}\label{eq:Jacelem37}
P^+=-\frac{(-1)^n}{2^n\,w(x)}\frac{e^{-\kappa\phi(z_+)}}{2\pi i}\int_{-\infty}^\infty e^{-\frac12\kappa s^2} f^+(s)\,ds,
\end{equation}
where
\begin{equation}\label{eq:Jacelem38}
f^+(s)=\frac{w(z)}{(z-x)}\left(\frac{z-x}{1-z^2}\right)^\gamma\frac{dz}{ds}.
\end{equation}
By expanding $\dsp{f^+(s)=f^+(0)\sum_{k=0}^\infty c_k^+s^k}$, the following expansion is obtained:
\begin{equation}\label{eq:Jacelem39}
P^+\sim A^+\sum_{m=0}^\infty \frac{ 2^m \left(\frac12\right)_mc_{2m}^+}{\kappa^m},
\quad A^+=-\frac{(-1)^n}{2^n\,w(x)}\frac{e^{-\kappa\phi(z_+)}}{2\pi i}f^+(0),
\end{equation}

For $A^+$ we evaluate
%%%
\begin{equation}\label{eq:Jacelem40}
\begin{array}{@{}r@{\;}c@{\;}l@{}}
&&\dsp{\frac{w(z_+)}{w(x)}=\frac{e^{\frac12i(\alpha\theta+\beta\theta-\pi\alpha)}}{\sin^{\alpha}\frac12\theta\,\cos^{\beta}\frac12\theta},}\quad\quad \dsp{ \frac{1-z_+^2}{z_+-x}=-2z_+,}\\[8pt]
&&\dsp{ \left. \frac{dz}{ds}\right\vert_{s=0}=\sqrt{\sin\theta}\,e^{\frac12i\theta-\frac14\pi i},}
\end{array}
\end{equation}
and this gives
\begin{equation}\label{eq:Jacelem41}
A^+=\frac{e^{i\chi}}{{2\sqrt{\pi\kappa}\, \sin^{\alpha+\frac12}\frac12\theta\,\cos^{\beta+\frac12}\frac12\theta}},
\end{equation}
with $\chi$ defined in \eqref{eq:Jacelem04}.

The contribution from the saddle point $z_-$ is the complex conjugate of $P^+$ (assuming that $\alpha$ and $\beta$ are real), and by splitting  up the coefficients of the expansion in \eqref{eq:Jacelem39} in real and imaginary parts, we obtain \eqref{eq:Jacelem03} and \eqref{eq:Jacelem05}.

\begin{remark}\label{rem:rem02}
The starting  integrand in \eqref{eq:Intro04} has a pole at $z=x$, while the one of \eqref{eq:Jacelem33} shows an algebraic singularity at $z=x$ and $\phi(z)$ defined in \eqref{eq:Jacelem34} has a logarithmic singularity at this point. To handle this from the viewpoint of multi-valued functions, we  can introduce a branch cut for the functions involved from $z=x$ to the left, assuming that the phase of $z-x$ is zero when $z>x$, equals $-\pi$ when $z$ approaches $-1$ on the lower part of the saddle point contour in \eqref{eq:Jacelem33}, and $+\pi$ on the upper side. Because the saddle points $e^{\pm i\theta}$ stay off the interval $(-1,1)$, we do not need to consider function values on the branch cuts for the asymptotic analysis.
\end{remark}

\section{An expansion in terms of Bessel functions}\label{sec:JacpolBess}

The second expansion we consider is expressed 
in terms of Bessel functions. The expansion is derived in \cite{Frenzen:1985:AUA} and 
in \S\ref{sec:BesAm} we give more details on the coefficients of this expansion. This expansion
will be used for the extreme zeros and, combined with the elementary expansion for the inner nodes, it 
will lead to a non-iterative double precision algorithm for Gauss-Jacobi quadrature. The expansion reads
\begin{equation}\label{eq:JacBes01}
P_n^{(\alpha,\beta)}(\cos\theta)\sim\frac{G_\kappa(\alpha,\beta)}{\sin^{\alpha}\frac12\theta\,
\cos^{\beta}\frac12\theta}\sqrt{\frac{\theta}{\sin\theta}}
\sum_{m=0}^\infty A_m(\theta) \frac{J_{\alpha+m}(\kappa\theta)}{\kappa^{m}},
\end{equation}
where  $G_\kappa(\alpha,\beta)$ is defined in \eqref{eq:Jacelem07},  and
\begin{equation}\label{eq:JacBes02}
\kappa=n+\tfrac12(\alpha+\beta+1).
\end{equation}
This holds uniformly valid with respect to $\theta\in[0,\pi-\delta]$, with $\delta$ an arbitrarily small positive number. 
The coefficients $A_m(\theta)$ are analytic functions for $0\le\theta<\pi$, and the first ones are $A_0(\theta)=1$ and
\begin{equation}\label{eq:JacBes03}
A_1(\theta)=\frac{\left(4\alpha^2-1\right)(\sin\theta-\theta\cos\theta)+2\theta\left(\alpha^2-\beta^2\right)(\cos\theta-1)}{8\theta\sin\theta}.
\end{equation}

We convert the expansion into the representation (see also \cite{Baratella:1988:TBF} and  \cite[Eqn.~18.15.6]{Koornwinder:2010:OPS})
\begin{equation}\label{eq:JacBes04}
\begin{array}{@{}r@{\;}c@{\;}l@{}}
P_n^{(\alpha,\beta)}(\cos\theta)&=&\dsp{\frac{G_\kappa(\alpha,\beta)}{\sin^{\alpha}\frac12\theta\,
\cos^{\beta}\frac12\theta}\sqrt{\frac{\theta}{\sin\theta}}\,W(\theta)}, \\[8pt]
W(\theta)&=&\dsp{J_{\alpha}(\kappa\theta)\,S(\theta)+\frac{1}{\kappa}J_{\alpha+1}(\kappa\theta)\,T(\theta),}
\end{array}
\end{equation}
where $S(\theta)$ and $T(\theta)$ have the expansions
\begin{equation}\label{eq:JacBes05}
S(\theta)\sim \sum_{m=0}^\infty\frac{S_m(\theta)}{\kappa^{2m}},\quad
T(\theta)\sim \sum_{m=0}^\infty\frac{T_m(\theta)}{\kappa^{2m}},\quad \kappa\to\infty,
\end{equation}
with $S_0(\theta)=A_0(\theta)=1$,  $T_0(\theta)=A_1(\theta)$, and for $m=1,2,3,\ldots$
\begin{equation}\label{eq:JacBes06}
\begin{array}{@{}r@{\;}c@{\;}l@{}}
S_m(\theta)&=&\dsp{-\frac{1}{\theta^{m-1}}\sum_{j=0}^{m-1} \binomial{m-1}{j}A_{j+m+1}(\theta)(-\theta)^j2^{m-1-j}(\alpha+2+j)_{m-j-1}},\\[8pt]
T_m(\theta)&=&\dsp{\frac{1}{\theta^{m}}\sum_{j=0}^{m} \binomial{m}{j}A_{j+m+1}(\theta)(-\theta)^j2^{m-j}(\alpha+1+j)_{m-j}}.  
\end{array}
\end{equation}
The first few terms are also given in \cite{Baratella:1988:TBF}, where the expansion is derived by using the differential equation of the Jacobi polynomials.

To compute the  coefficients $A_m(\theta)$ for small values of $\theta$ in a stable way we need expansions. We can write
\begin{equation}\label{eq:JacBes07}
A_m(\theta)=\chi^m \theta^{m}\sum_{j=0}^\infty A _{jm}\theta^{2j},\quad \chi=\frac{\theta}{\sin\theta},
\end{equation}
where the series represent entire functions of $\theta$. The first few $A_{jm}$ are
\begin{equation}\label{eq:JacBes08}
\begin{array}{@{}r@{\;}c@{\;}l@{}}
A_{0,1}&=& \dsp{\tfrac{1}{24}\left(\alpha^2+3\beta^2-1\right), }\\[8pt]
A_{1,1}&=& \dsp{\tfrac{1}{480}\left(-3\alpha^2-5\beta^2+2\right), }\\[8pt]
A_{0,2}&=& \dsp{\tfrac{1}{5760}\left(-16\alpha-14\alpha^2-90\beta^2+5\alpha^4+4\alpha^3+45\beta^4+30\beta^2\alpha^2+60\beta^2\alpha+21\right). }
\end{array}
 \end{equation}

The coefficients $S_m(\theta)$ and $T_m(\theta)$ can be expanded in the form
\begin{equation}\label{eq:JacBes09}
\begin{array}{@{}r@{\;}c@{\;}l@{}}
S_m(\theta)&=&\dsp{\theta^2\chi^{2m} \sum_{j=0}^\infty S _{jm}\theta^{2j},\quad m\ge1,}\\[8pt]
T_m(\theta)&=&\dsp{\theta\chi^{2m+1} \sum_{j=0}^\infty T _{jm}\theta^{2j},\quad m\ge0,}
\end{array}
 \end{equation}
in which the series represent entire functions of $\theta$.

\begin{remark}\label{rem:rem03}
Asymptotic expansions of Jacobi polynomials in terms of modified Bessel functions $I_\nu(z)$  have been derived in \cite{Elliott:1971:UAE}, where $P_n^{(\alpha,\beta)}(z)$ and an associated function (with expansion in terms of $K_\nu(z)$) are considered for $z\notin[-1,1]$. The approach is based on Frank Olver's methods  by using the differential equation. However, Elliott's expansion of $P_n^{(\alpha,\beta)}(z)$ is also valid for $z$ in a neighborhood of $z=1$, and it coincides with our expansion when performing the proper transformations.
%\eoremark
\end{remark}

\subsection{Expansions of a derivative}\label{sec:JacBesder}
For computing the Gauss weights it is convenient to have an expansion of the derivative of the function $U(\theta)=\sqrt{\theta}\,W(\theta)$, with $W(\theta)$ defined in \eqref{eq:JacBes04}.
We have
\begin{equation}\label{eq:JacBes12}
\frac{d}{d\theta}U(\theta)=-\kappa\sqrt{\theta}\left(J_{\alpha+1}(\kappa\theta)\,Y(\theta)-\frac{1}{2\theta\kappa}J_{\alpha}(\kappa\theta)\,Z(\theta)\right),
\end{equation}
with expansions
\begin{equation}\label{eq:JacBes13}
Y(\theta)\sim \sum_{m=0}^\infty\frac{Y_m(\theta)}{\kappa^{2m}},\quad
Z(\theta)\sim \sum_{m=0}^\infty\frac{Z_m(\theta)}{\kappa^{2m}},\quad \kappa\to\infty,
\end{equation}
where $Y_0(\theta)=S_0(\theta)=A_0(\theta)=1$,  $Z_0(\theta)=(2\alpha+1)+2\theta A_1(\theta)$, and for $m=1,2,3,\ldots$
\begin{equation}\label{eq:JacBes14}
\begin{array}{@{}r@{\;}c@{\;}l@{}}
Y_m(\theta)&=&\dsp{S_m(\theta)+\frac{2\alpha+1}{2\theta}T_{m-1}(\theta)-\frac{d}{d\theta}T_{m-1}(\theta)},\\[8pt]
Z_m(\theta)&=&\dsp{(2\alpha+1)S_m(\theta)+2\theta T_m(\theta)+2\theta\frac{d}{d\theta}S_m(\theta)}.  
\end{array}
\end{equation}

For small values of $\theta$ we need expansions, similar as in \eqref{eq:JacBes07}. We have
\begin{equation}\label{eq:JacBes15}
\begin{array}{@{}r@{\;}c@{\;}l@{}}
Y_m(\theta)&=&\dsp{\chi^{2m} \sum_{j=0}^\infty Y _{jm}\theta^{2j},\quad m\ge1,}\\[8pt]
Z_m(\theta)&=&\dsp{\theta^2\chi^{2m+1} \sum_{j=0}^\infty Z _{jm}\theta^{2j},\quad m\ge0,}
\end{array}
 \end{equation}
in which the series represent entire functions of $\theta$.

\subsection{Expansions of the zeros}\label{sec:JacBeszer}

For obtaining accurate approximations of the zeros for large degree Jacobi polynomials we can 
use the Bessel-type expansion given earlier. We give asymptotic expansions that can be used for all positive zeros. However, 
the most interesting region of application of Bessel-type expansions is close to the endpoints of the interval
$[-1,1]$, because for the rest of the interval the previous elementary expansions are accurate enough 
and they are simpler to handle (and therefore more efficient). 
We explain how to compute explicitly enough coefficients
of the expansions so that the zeros are computed accurately, without the need to use iterative methods to refine the results.

We write the zeros $x_{n+1-m}$ (with $x_1<x_2<\cdots<x_n$)
of $P_n^{(\alpha,\beta)}(x)$ in terms of the zeros $j_m$ of the Bessel function $J_\alpha(x)$. The zero $x_n$ 
corresponds to $j_1$, $x_{n-1}$ to $j_2$, and so on. Because the representation in \eqref{eq:JacBes04} cannot be used as $x\to-1$, we consider only the positive zeros. 
For the other zeros we can use the symmetry relation \eqref{eq:Intro01}.

A zero of $P_n^{(\alpha,\beta)}(x)$ is a zero of $W(\theta)$ defined in \eqref{eq:JacBes04}, and we approximate a zero $x_{n+1-m}$ with corresponding  $\theta$ value following from $x_{n+1-m}=\cos\theta$. We write $\theta$ in the form
\begin{equation}\label{eq:JacBes16}
\theta=\theta_0+\eps,\quad \theta_0=j_m/\kappa.
\end{equation}
We assume for $\eps$ an expansion in the form
\begin{equation}\label{eq:JacBes17}
\eps\sim\frac{\theta_1}{\kappa^2}+\frac{\theta_2}{\kappa^4}+\frac{\theta_3}{\kappa^6}+\ldots.
\end{equation}
By expanding $W(\theta)$ at $\theta_0$ we have upon expanding
\begin{equation}\label{eq:JacBes18}
W(\theta)=W(\theta_0+\eps)=W(\theta_0)+\frac{\eps}{1!}W^\prime(\theta_0)+ \frac{\eps^2}{2!}W^{\prime\prime}(\theta_0)+\ldots = 0.
\end{equation}
Using the representation of $W(\theta)$ given in \eqref{eq:JacBes04},  substituting the expansion of $\eps$ and those of $S(\theta)$ and $T(\theta)$ given in \eqref{eq:JacBes05}, and comparing equal powers of $\kappa$, we can obtain the coefficients $\theta_j$ of \eqref{eq:JacBes17}. 

The first coefficient is $\theta_1=T_0(\theta_0)=A_1(\theta_0)$ (see \eqref{eq:JacBes09}), and the second one follows from
\begin{equation}\label{eq:JacBes19}
6\theta\theta_2=6\theta T_1(\theta)+6\theta(T_0^\prime(\theta)-S_1(\theta)) T_0(\theta)-3(1+2\alpha)T_0(\theta)^2-2\theta T_0(\theta)^3,
\end{equation}
evaluated at $\theta=\theta_0$.

The coefficients can be expanded in the form
\begin{equation}\label{eq:JacBes20}
\theta_m=\theta\chi^{2m-1} \sum_{j=0}^\infty t_{jm}\theta^{2j},\quad m\ge1,
\end{equation}
evaluated at $\theta_0$, where the series represent entire functions of $\theta$. The first few $t_{jm}$ are 
\begin{equation}\label{eq:JacBes21}
\begin{array}{@{}r@{\;}c@{\;}l@{}}
t_{0,2} &=& \frac12(2T_{0, 1}+T_{0, 0}^2-2\alpha T_{0, 0}^2),\\[8pt]
t_{1,2}&=&\frac{1}{12}(-4T_{0, 0}^3+12T_{1, 1}-12T_{0, 0}S_{0, 1}-24\alpha T_{0, 0}T_{1, 0}\ +\\[8pt]
&&\quad\quad
36T_{0, 0}T_{1, 0}+3T_{0, 0}^2+2\alpha T_{0, 0}^2),\\[8pt]
t_{0,3}&=&\frac{1}{12}(12T_{0, 2}+12T_{0, 0}T_{0, 1}+2T_{0, 0}^3\ +\\[8pt]
&&\quad\quad
16\alpha^2T_{0, 0}^3-12\alpha T_{0, 0}^3-24\alpha T_{0, 0}T_{0, 1}).
\end{array}
 \end{equation}

The numerical steps to obtain the zeros are as explained at the end of \S\ref{sec:Jaczerelem}. 
For example, with $n=100$, $\alpha=\frac13$, $\beta=\frac14$, we have for the largest zero $x_{100}=0.9995853721163790$ (computed by Maple), $\theta_0=0.02879787927325625$, and an approximation of $x_{100}$ given by $\cos\theta_0=0.9995853697308934$, with relative error 
$2.39\times10^{-9}$. Compute $\theta_1=-0.8416536425087086\times10^{-3}$, and we have $\theta\sim\theta_0+\theta_1/\kappa^2$, which gives $x_{100}\sim\cos\theta=0.9995853721164185$, with relative error $3.96\times10^{-14}$. These computations were done in Maple with Digits=16.

\subsection{Details on computing Bessel functions near a zero}\label{sec:BesBes}
The numerical evaluation of the Bessel  function $J_{\alpha}(\kappa\theta)$ that occurs in $W(\theta)$ defined in \eqref{eq:JacBes04} and in the derivative of $U(\theta)$ in \eqref{eq:JacBes12}, may become unstable when $\kappa\theta$ is near a zero of $J_{\alpha}(z)$. This problem shows up when computing the Gauss weights and when we use the standard software for the Bessel functions. A related problem is discussed in \S\ref{sec:Jaccos} for the elementary case.

To handle this, we can use the following expansion (see \cite[Eqn.~10.23.1]{Olver:2010:Bessel})
\begin{equation}\label{eq:JacBes22}
J_\alpha(\lambda z)=\lambda^\alpha\sum_{m=0}^\infty \frac{w^m}{m!}J_{\alpha+m}(z),\quad w=\tfrac12z\left(1-\lambda^2\right),
\end{equation}
with $z=u$ and $\lambda=1+h/u$. This gives
\begin{equation}\label{eq:JacBes23}
J_\alpha(u+h)=\lambda^\alpha\sum_{m=0}^\infty \frac{w^m}{m!}J_{\alpha+m}(u),\quad w=-\frac{h(2u+h)}{2u}.
\end{equation}
When $h$ is small the series works as a Taylor expansion, because $w=\bigO(h)$. When $u$ is a zero of $J_\alpha(u)$, the first term of the series vanishes.

This happens when we take $u+h=\kappa\theta$ with $u=\kappa\theta_0=j_k$ (a zero of $J_\alpha(u)$, see \eqref{eq:JacBes16}) and $h=\kappa\eps$. The mentioned computational problems arise for largest zeros $x_k$ (with small $k$). In that case, $u=\bigO(1)$ and $h=\bigO(1/\kappa)$, and the series in \eqref{eq:JacBes23} converges quite fast. 

The terms of the series can be generated by using the recurrence relation of the Bessel function
\begin{equation}\label{eq:JacBes24}
J_{\alpha+1}(u)=\frac{2\alpha}{u}J_{\alpha}(u)-J_{\alpha-1}(u),
\end{equation}
which gives for $m=1,2,3,\ldots$
\begin{equation}\label{eq:JacBes25}
m(m+1)f_{m+1}=\frac{2m(\alpha+m)w}{u}f_m-w^2f_{m-1}, \quad f_m=\frac{w^m}{m!}J_{\alpha+m}(u),
\end{equation}
with starting value $f_0=0$ (when $u$ is a zero of $J_{\alpha}(u)$) and $f_1=wJ_{\alpha+1}(u)$.

In Table~\ref{tab:tab01} we show the relative errors in the computation of the Bessel function $J_{\alpha}(u+h)$ for $\alpha=\frac14$, $u=j_5=15.32\cdots$ and several small values of $h$ by using the standard Maple code for this Bessel function and the series expansion in \eqref{eq:JacBes23}, both with Digits = 16. The value $m$ indicates the number of terms used in the expansion. The errors follow from comparing the results with Bessel function values obtained in Maple with Digits = 64. We observe a good performance of the series expansions, the standard algorithm gives poorer results as $h$ becomes smaller.

\begin{remark}\label{rem:rem04}
We know that the forward recursion of the Bessel functions in \eqref{eq:JacBes25} may be unstable. However, because we start in the domain where the zeros are ($u>\alpha$), the recursion of the first terms will be stable, and for the later terms we can accept 
less accurate values. Because of the fast convergence of the series for the values of $u$ and $h$ to compute the Gauss weights, it is not needed to use a backward recursion scheme for the Bessel functions.
%\eoremark
\end{remark}

\renewcommand{\arraystretch}{1.2}
\begin{table}%[H]
\caption{
Relative errors in the computation of the Bessel function $J_{\alpha}(u+h)$ for $\alpha=\frac14$, $u=j_5=15.32\cdots$ and several small values of $h$ by using the standard Maple code for this Bessel function with Digits = 16 and the series expansion in \eqref{eq:JacBes23}; $m$ indicates the number of terms used in the expansion. 
\label{tab:tab01}}
$$
\begin{array}{cccc}
h & m& {\rm standard\  algorithm} \ & {\rm series\  expansion}\  \eqref{eq:JacBes23}  \\
\hline
10^{-1}\quad &\quad 9 \quad & \quad0.26\times 10^{-13}\quad & \quad0.26\times 10^{-15}\\
10^{-2}\quad &\quad 5 \quad & \quad0.26\times 10^{-12}\quad & \quad0.68\times 10^{-16}\\
10^{-3}\quad &\quad 4 \quad & \quad0.26\times 10^{-11}\quad & \quad0.17\times 10^{-15}\\
10^{-4}\quad &\quad 3 \quad & \quad0.26\times 10^{-10}\quad & \quad0.55\times 10^{-15}\\
10^{-5}\quad &\quad 2 \quad & \quad0.26\times 10^{-09}\quad & \quad0.13\times 10^{-15}\\
\hline
\end{array}
$$
\end{table}
\renewcommand{\arraystretch}{1.0}

\subsection{Details on the coefficients $A_m(\theta)$}\label{sec:BesAm}

We verify the expansion in \eqref{eq:JacBes01} and give details on the coefficients $A_m(\theta)$. We introduce
\begin{equation}\label{eq:JacBes26}
b_j(\theta)=\frac{(-1)^j\left(\frac12\alpha+\frac12\beta\right)_j\left(\frac12\alpha-\frac12\beta\right)_j}
{j!\,\left(\alpha+\frac12\right)_j\left(1+\cos\theta\right)^j},\quad j=0,1,2,\ldots, 
\end{equation}
and for $m=0,1,2,\ldots$ functions $\psi_{j,m}(\theta)$ defined by the generating functions
\begin{equation}\label{eq:JacBes27}
\left(\frac{2\theta(\cos\phi-\cos\theta)}{\sin\theta\left(\theta^2-\phi^2\right)}\right)^{m+\alpha-\frac12}=
\sum_{j=0}^\infty \psi_{j,m}(\theta)\left(\theta^2-\phi^2\right)^{j}. 
\end{equation}
Then, the coefficients $A_m(\theta)$ are  given by
\begin{equation}\label{eq:JacBes28}
A_m(\theta)=(2\theta)^{m}\left(\alpha+\tfrac12\right)_{m}\sum_{j=0}^{m}b_j(\theta)\left(\frac{\sin\theta}{2\theta}\right)^j\psi_{m-j,j}(\theta).
\end{equation}

The starting point in \cite{Frenzen:1985:AUA} is the integral representation (due to George Gasper \cite{Gasper:1975:FDM}) in terms of the Gauss hypergeometric function
\begin{equation}\label{eq:JacBes29}
\begin{array}{@{}r@{\;}c@{\;}l@{}}
&&P_n^{(\alpha,\beta)}(\cos\theta)=Q_n^{(\alpha,\beta)}(\theta)\times\\[8pt]
&&\quad
\dsp{\int _0^\theta
\frac{\cos(\kappa\phi)}{(\cos \phi-\cos\theta)^{\frac12-\alpha}}
\FG{\frac12(\alpha+\beta)}{\frac12(\alpha-\beta)}{\alpha+\frac12}
{\frac{\cos\theta-\cos\phi}{1+\cos\theta}}\,d\phi,}
\end{array}
\end{equation}
where $\kappa$ is defined in \eqref{eq:JacBes02} and
\begin{equation}\label{eq:JacBes30}
Q_n^{(\alpha,\beta)}(\theta)=\frac{2^{\frac12-\alpha}\Gamma(n+\alpha+1)}{\sqrt\pi\,n!\,\Gamma\left(\alpha+\frac12\right)
\sin^{2\alpha}\left(\frac12\theta\right) \cos^{\alpha+\beta}\left(\frac12\theta\right)}.
\end{equation}

By expanding the ${}_2F_1$-function, it follows that
\begin{equation}\label{eq:JacBes31}
\begin{array}{@{}r@{\;}c@{\;}l@{}}
P_n^{(\alpha,\beta)}(\cos\theta)&=&Q_n^{(\alpha,\beta)}(\theta)\times\\[8pt]
&&\dsp{\sum_{m=0}^\infty b_m\int _0^\theta
\cos(\kappa\phi)(\cos \phi-\cos\theta)^{m+\alpha-\frac12}}\,d\phi,
\end{array}
\end{equation}
where $b_m(\theta)$ is defined in \eqref{eq:JacBes26}

%%5
Next, the expansion in \eqref{eq:JacBes27} is used and the integral representation\footnote{A proof easily follows  by expanding the cosine function in its power series.}
\begin{equation}\label{eq:JacBes32}
\int_0^\theta\left(\theta^2-\phi^2\right)^{\sigma-\frac12}\cos(\kappa\phi)\,d\phi=2^{\sigma-1}\sqrt\pi\,\Gamma\left(\sigma+\tfrac12\right)
\theta^\sigma\kappa^{-\sigma}J_\sigma(\kappa \theta),
\end{equation}
with $\sigma=j+m+\alpha$. This gives the double sum
\begin{equation}\label{eq:JacBes33}
\begin{array}{@{}r@{\;}c@{\;}l@{}}
P_n^{(\alpha,\beta)}(\cos\theta)
&=&\dsp{Q_n^{(\alpha,\beta)}(\theta)\sqrt{\frac{\pi}{2}}\,\left(\frac{\sin\theta}{\theta}\right)^{\alpha-\frac12}\left(\frac{\theta}{\kappa}\right)^{\alpha}}\times\\[8pt]
&&\dsp{\sum_{m=0}^\infty p_m\sum_{j=0}^\infty q_{j+m}\psi_{j,m}(\theta)}\\[8pt]
&=&\dsp{Q_n^{(\alpha,\beta)}(\theta)\sqrt{\frac{\pi}{2}}\,\left(\frac{\sin\theta}{\theta}\right)^{\alpha-\frac12}\left(\frac{\theta}{\kappa}\right)^{\alpha}}\times\\[8pt]
&&\dsp{\sum_{m=0}^\infty q_m\sum_{j=0}^m p_j \psi_{m-j,j}(\theta)},
\end{array}
\end{equation}
where $\psi_{j,m}(\theta)$ are the coefficients  in \eqref{eq:JacBes27}, and
\begin{equation}\label{eq:JacBes34}
p_m=b_m(\theta)\left(\frac{\sin\theta}{2\theta}\right)^m, \quad
q_{m}=\left(\frac{2\theta}{\kappa}\right)^{m}\Gamma\left(m+\alpha+\tfrac12\right)J_{m+\alpha}(\kappa\theta).
\end{equation}

This verifies the expansion in \eqref{eq:JacBes01} with the $A_m(\theta)$ defined in \eqref{eq:JacBes28}.

We give a few details from \cite{Frenzen:1985:AUA} on evaluating the functions $\psi_{j,m}(\theta)$. First consider the expansion (see \cite[p.~140 (3)]{Watson:1944:TTB})
\begin{equation}\label{eq:JacBes35}
\frac{2\theta(\cos\phi-\cos\theta)}{\sin\theta\left(\theta^2-\phi^2\right)}=
\sum_{j=0}^\infty \chi_{j}(\theta)\left(\theta^2-\phi^2\right)^{j},
\end{equation}
where the $\chi_{j}(\theta)$ are available in the form of Bessel functions of fractional order:
\begin{equation}\label{eq:JacBes36}
\chi_{j}(\theta)=\frac{1}{(j+1)!\,(2\theta)^j}\frac{J_{j+\frac12}(\theta)}{J_{\frac12}(\theta)},\quad j=0,1,2,,\ldots.
\end{equation}
These functions follow from a simple recursion, and then the $\psi_{j,m}(\theta)$ follow from
\begin{equation}\label{eq:JacBes37}
\left(\sum_{j=0}^\infty \chi_{j}(\theta)w^{j}\right)^{\mu}=\sum_{j=0}^\infty \psi_{j,m}(\theta)w^{j},\quad w= \theta^2-\phi^2, \quad\mu=m+\alpha-\tfrac12.
\end{equation}
The first two are
\begin{equation}\label{eq:JacBes38}
\psi_{0,m}(\theta)=1,\quad \psi_{1,m}(\theta)=\tfrac{1}{4}\mu\frac{1-\theta\cot\theta}{\sin^2\theta}.
\end{equation}

\section{Numerical performance of the expansions for computing the nodes and weights of the G-J quadrature}\label{sec:numerical}

Next, we discuss in more detail the performance of the expansions given in \S\ref{sec:Jacpelem} and 
 \S\ref{sec:JacpolBess}
for computing the nodes
and weights of the Gauss-Jacobi quadrature. 

In terms of the derivatives of the orthogonal polynomials, the weights of the Gauss-Jacobi quadrature are given by

\begin{equation}
\label{eq:numer1}
\begin{array}{l}
w_i  =
\Frac{M_{n,\alpha,\beta}}{ (1-x_i^2) [P_n ^{(\alpha ,\beta) \prime}(x_i)]^2} = \Frac{M_{n,\alpha,\beta}}{\left[\Frac{d}{d\theta}P_n ^{(\alpha ,\beta) }(\cos \theta_i)\right]^2},\\
\\
M_{n,\alpha,\beta}=2^{\alpha+\beta+1}\Frac{\Gamma (n+\alpha+1)\Gamma (n+\beta+1)}{n! \Gamma (n+\alpha+\beta+1 )},
\end{array}
\end{equation}
where $x_i=\cos\theta_i$.

As we did in \cite{Gil:2018:GHL} for the Gauss--Hermite and Gauss--Laguerre quadratures, it is convenient 
to introduce scaled weights: 

\begin{equation}
\label{eq:numer2}
\omega_i=\dot{u}(\theta_i)^{-2}
\end{equation}
where the dot indicates the derivative with respect to $\theta$ and
 
\begin{equation}
\label{eq:numer3}
u(\theta)=M^{-\frac12}_{n,\alpha,\beta}\left(\sin \frac12\theta\right)^{\alpha +\frac12}
\left(\cos \frac12\theta\right )^{\beta +\frac12} P_n^{(\alpha,\beta)}(\cos\theta).
\end{equation}

The weights are related with the scaled weights by
\begin{equation}
\label{relat}
w_i =\left(\sin \frac12\theta\right)^{2\alpha +1}
\left(\cos \frac12\theta\right )^{2\beta +1} \omega_i .
\end{equation}

The advantage of computing scaled weights is that, similarly as described in \cite{Gil:2018:GHL},
scaled weights do not underflow/overflow for large parameters and, additionally, they are well-conditioned
as a function of the roots $\theta_i$. Indeed, the scaled weights are $\omega_i=W(\theta_i)$ with
$W(\theta )=\dot{u}(\theta)^{-2}$, but $\dot{W}(\theta_i)=0$ because $\ddot{u}(\theta_i)=0$. In addition,
as $n\rightarrow +\infty$, the scaled weights are essentially constant and 
the main dependence on the nodes of the unscaled weights $w_i$ is given by the elementary sine and 
cosine factor of (\ref{relat}). This is related to the
circle theorem for Gaussian quadratures in $[-1,1]$ 
(see \cite{Davis:1961:SGT,Gautschi:2006:TCT}), 
which states that, as $n\rightarrow +\infty$,
$$
\Frac{n w_i}{\pi w(x_i)}\sim\sqrt{1-x_i^2},
$$
where in the Gauss-Jacobi case $w(x)=(1-x)^\alpha (1+x)^\beta$. 

When considering the asymptotic expansion for Jacobi polynomials in terms of elementary functions,
the function $u(\theta)$ can be written in terms of the function $W(\theta)$ given in (\ref{eq:Jacelem19})

\begin{equation}
\label{eq:numer4}
u(\theta)=
\Frac{M^{-\frac12}_{n,\alpha,\beta}G(\alpha,\,\beta)}{\sqrt{\pi \kappa} }W(\theta),
\end{equation}
and for the computation of $\dot{u}(\theta)$ we use the expansion (\ref{eq:Jacelem20}).

When considering the asymptotic expansion in terms of Bessel functions, the function
$u(\theta)$ is given by

\begin{equation}
\label{eq:numer5}
u(\theta)=
\Frac{M^{-\frac12}_{n,\alpha,\beta}\Gamma(n+\alpha+1)}{n!\kappa^{\alpha}\sqrt{2}}U(\theta),
\end{equation}
where $U(\theta)$ is given in Section \ref{sec:JacBesder}.
In this case, the expansion \eqref{eq:JacBes12} is used for computing 
the derivative of $u(\theta)$.

\begin{figure}
\epsfxsize=15cm \epsfbox{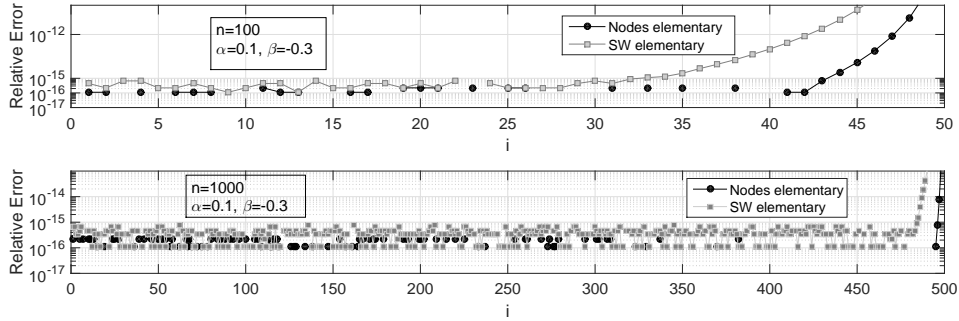}
\caption{
\label{fig:fig01} Performance of the asymptotic expansions in terms of
elementary functions for computing the nodes
and scaled weights of $P^{(\alpha,\beta)}_n(x)$ for
 $\alpha=0.1$, $\beta=-0.3$ and $n=100,\,1000$.}
\end{figure}

\begin{figure}
\epsfxsize=15cm \epsfbox{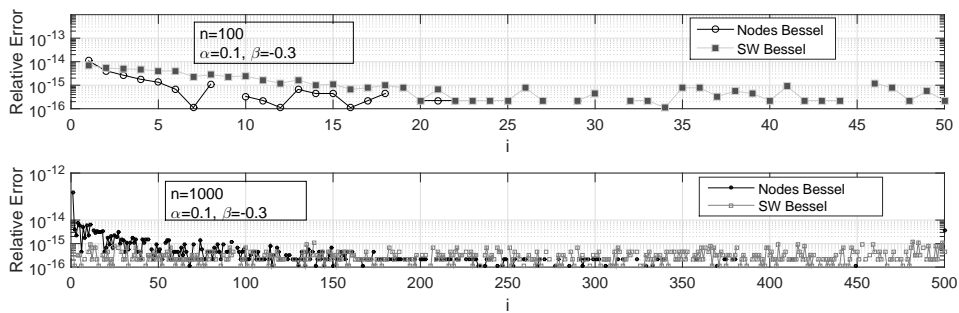}
\caption{
\label{fig:fig02} Performance of the asymptotic expansions
in terms of Bessel functions  for computing the nodes
and scaled weights of $P^{(\alpha,\beta)}_n(x)$ for
 $\alpha=0.1$, $\beta=-0.3$ and $n=100,\,1000$.}
\end{figure}

\begin{figure}
\epsfxsize=15cm \epsfbox{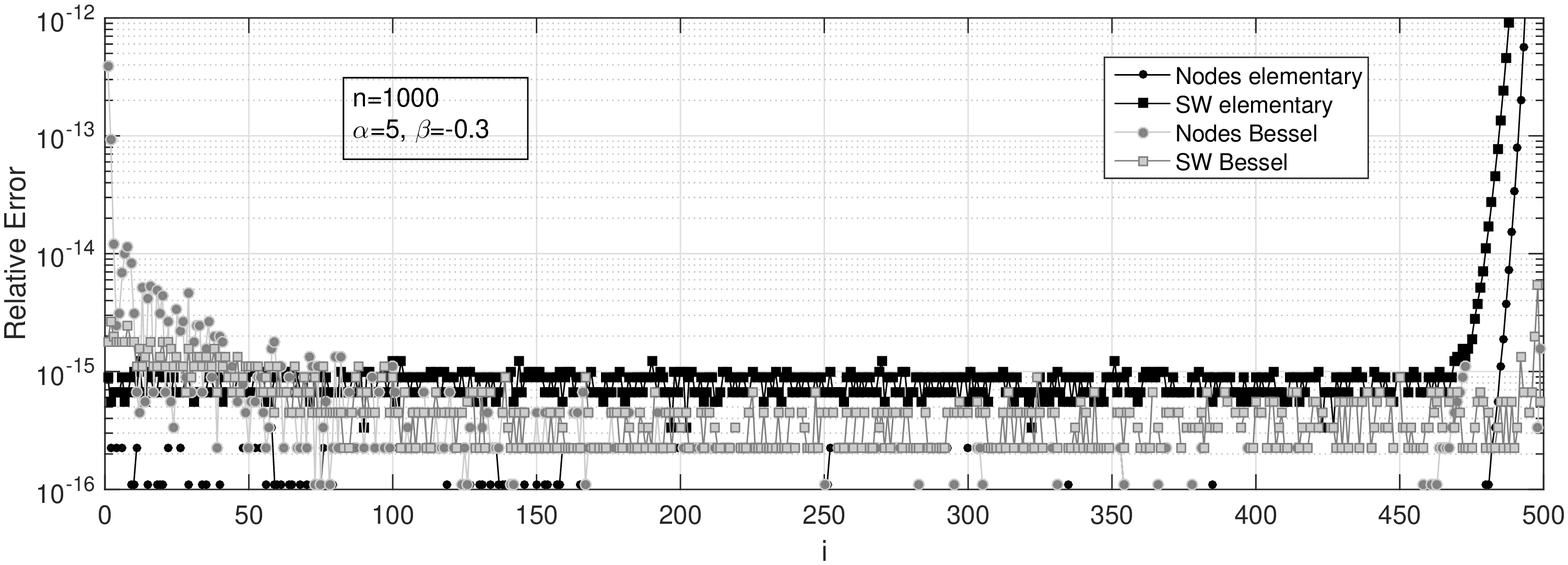}
\caption{
\label{fig:fig03} Performance of the asymptotic expansions for computing the nodes
and scaled weights of $P^{(\alpha,\beta)}_n(x)$ for
 $\alpha=5$, $\beta=-0.3$ and $n=1000$.
}
\end{figure}

\begin{figure}
\epsfxsize=15cm \epsfbox{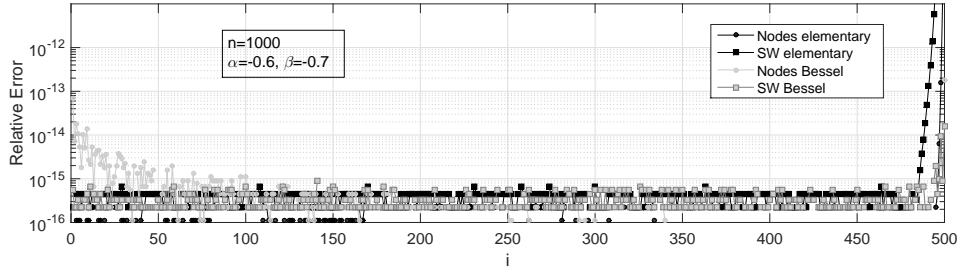}
\caption{
\label{fig:fig04} Performance of the asymptotic expansions for computing the nodes
and scaled weights of $P^{(\alpha,\beta)}_n(x)$ for
 $\alpha=-0.6$, $\beta=-0.7$ and $n=1000$.
}
\end{figure}

Examples of the performance of the  asymptotic expansions for the evaluation of the nodes
and scaled weights of the Gauss-Jacobi quadrature are given
 in Figures \ref{fig:fig01}, \ref{fig:fig02}, \ref{fig:fig03} and \ref{fig:fig04}. We concentrate on the positive zeros;
for the negative zeros we can use the relation (\ref{eq:Intro01}). 
In our tests, we use finite
precision implementations (double precision accuracy) of the expansions
\footnote{A file with explicit expressions of some of the coefficients of the expansions
used in our calculations is available at http://personales.unican.es/gila/CoefsExpansions.txt; see also 
http://personales.unican.es/segurajj/gaussian.html}.
We compare the 
approximations to the nodes obtained with the asymptotic expansions 
against the results of a Maple implementation (with a large number of digits)
of an iterative algorithm which uses the global fixed point method of 
\cite{Segura:2010:RCO}. The Jacobi
polynomials used in this algorithm are computed using the intrinsic Maple function.
The scaled weights for testing have also been computed using Maple. 

 Figure \ref{fig:fig01} shows the performance of the expansion 
in terms of elementary functions for $P^{(\alpha,\beta)}_n(x)$ for
 $\alpha=0.1$, $\beta=-0.3$ and $n=100,\,1000$. 
For computing the scaled weights (SW in the figure) for $n=100$, five terms have been considered
in the expansions of the coefficients $M(x)$ and $N(x)$ in (\ref{eq:Jacelem21}).
For the nodes, we have used six terms in the expansion   (\ref{eq:Jacelem26}).
For $n=1000$, fewer terms were needed to obtain the accuracy shown in the figure.
The label $i$ in the abscissa represents the order of the zero (starting from $i=1$ for the smallest positive zero). 
 The points not shown in the plot correspond 
to values with all digits correct in double precision accuracy.
As can be seen,  for $n=100$
the use of the elementary expansion allows the computation of more than $3/4$ of the first zeros with 15-16
digits correct and it fails for the largest zeros, as expected. When $n=1000$, only the last few nodes
are not computed with double precision accuracy.
It is important to note that, using the idea given \S\ref{sec:Jaccos}
 for the computation of $\cos \theta$, there is no loss of accuracy for the first nodes and then double
precision accuracy is also possible in the relative error of these nodes, which is an advantage over
 the algorithms of \cite{Hale:2013:FAC}.

The asymptotic expansion for the scaled weights performs also well for the first half of the nodes ($n=100$)
with the number of terms used in the computation. For $n=1000$, only the last $15$ scaled weights are
computed with an accuracy worse than $10^{-15}$. 

The elementary expansions are fast and for $n=100$ the average CPU time per node 
(running gfortran on a Laptop with Intel Core i54310U 2.6 GHz processor under Windows) 
was $\approx 0.7 \mu s$, and the additional computation of  weights increased the
total CPU time by a factor 3. As the degree becomes larger, the time per node can be smaller because 
less terms are required, but only with this (pessimistic) estimation and taking into account we 
observe a clear speed-up with respect to the methods in \cite{Hale:2013:FAC}. This, is as expected, 
because for computing the weights we only need to evaluate 
one additional asymptotic, while the iterative methods must refine the nodes and two expansions (one for
the polynomial and another one for the derivative) must be computed. We therefore expect
that our computations will be around factor $2N$ faster, with $N$ the average number of iterations used in
\cite{Hale:2013:FAC}. 
%%%5   

For the computation of the nodes using the Bessel expansion, an algorithm for computing the zeros of Bessel
functions is needed. In the finite precision implementation of the expansion 
we use the algorithm given in \cite{Gil:2012:CZC}
for the first few zeros. For larger zeros, a very efficient method of computation is the use of the MacMahon's
expansion (see \cite[\S10.21(vi)]{Olver:2010:Bessel})

\begin{equation}
j_{\nu,m} \sim a-\frac{\mu-1}{8a}-\frac{4(\mu-1)(7\mu-31)}{3(8a)^{3}}%
-\frac{32(\mu-1)(83\mu^{2}-982\mu+3779)}{15(8a)^{5}}-\cdots,
\end{equation}
where $\mu=4\nu^2$, $a=(m+\nu/2-1/4)\pi$.

On the other hand, the $J$-Bessel functions needed to compute the weights, are evaluated using our own
algorithm for Bessel functions.  Two first examples of the performance of the Bessel expansion are shown in
Figure \ref{fig:fig02}. For comparison, the choice of parameters is the same as in Figure \ref{fig:fig01}.
For computing the scaled weights for $n=100$, four terms have been considered
in the expansions of the coefficients $Y(x)$ and $Z(x)$ in (\ref{eq:JacBes13}).
For the nodes, we have used three terms in the expansion   (\ref{eq:JacBes17}).
As can be seen, the expansion performs well even for the smallest nodes (and their corresponding scaled
weights), although higher accuracy is obtained for larger nodes as expected. 
Regarding the computation of scaled weights, it is important
to mention that the accurate evaluation of the two Bessel functions appearing in the expansion, seems
to be crucial in order to obtain more than $14$ digits of accuracy for a few of the points closest to
$x=1$ \footnote{In \cite{Hale:2013:FAC} it was reported that a loss of accuracy in the weights near the boundary took place,
and that it appeared to be directly related to the error in evaluating the Bessel functions at these points.}.
For computing the Bessel function $J_{\alpha}(\kappa\theta)$ the scheme given 
in \S\ref{sec:BesBes} is used in our calculations for the last points.

We have checked that the range validity of the two different expansions depend very little on the values of
$\alpha$ and $\beta$ ($-1<\alpha,\beta\le 5$) and that the main dependence is on the degree $n$. Because the
elementary expansion is more efficient, it should be preferred over the Bessel expansion when it gives
sufficiently accurate results. This is, as mentioned, the case of $3/4$ of all the nodes when $n=100$ 
and decreasing as $n$ increases. Once the nodes are computed, the weights are obtained, as explained before, 
by using the asymptotic expansions (\ref{eq:Jacelem20}) and (\ref{eq:JacBes12}); for this computation, it
seems convenient to use the elementary expansion (\ref{eq:Jacelem20}) for less nodes (half of the nodes for 
$n=100$), and more weights can be computed with this elementary expansion as $n$ becomes larger. 
%%%5

Two other examples of the computation of nodes and scaled weights 
with the expansions considered in this paper are shown in Figures \ref{fig:fig03} and \ref{fig:fig04}. 
In Figure \ref{fig:fig03}, the parameter $\alpha$ is chosen now larger than before
($\alpha=5$) in order to show the performance of the expansions for moderate values
of the parameters. On the other hand, Figure \ref{fig:fig04} illustrates
the performance when both parameters ($\alpha$ and $\beta$) are negative. 
 As can be seen, the performance of the expansions is very similar in all cases.

Regarding CPU times, the Bessel expansions, as expected, are slower than the elementary expansions. 
For $n=100$ the average CPU time per node  (running gfortran on an Intel Core i54310U 2.6 GHz processor under Windows) was $2\, \mu s$. The additional computation of  weights increased the total CPU time by a factor $5$. This expansion
is, as expected, slower than the elementary expansion. However, as mentioned earlier, the Bessel expansions can be
avoided for most of the nodes and the execution time of the algorithm combining both expansions is closer to the 
elementary expansion timing than to the Bessel case. The actual speed-up of the combined 
algorithm in the least favorable case ($n=100$) is larger than a factor $20$ with respect to the CPU time 
declared in Table 4.1 for the Matlab algorithm of \cite{Hale:2013:FAC}; 
however this comparison should be handled with care because 
it is always difficult to compare algorithms which are executed in different computers, using different 
programming languages, and probably with different optimization levels. A detailed discussion of the 
comparison of CPU times will be only possible when the algorithms are made available in a same platform. 
In any case, taking into account the complexities of the two methods, and as commented earlier, we can safely 
expect that our methods will be at least $2N$ times faster, 
with $N$ the average number of iterations used by the algorithm in \cite{Hale:2013:FAC}.

\begin{figure}
\epsfxsize=15cm \epsfbox{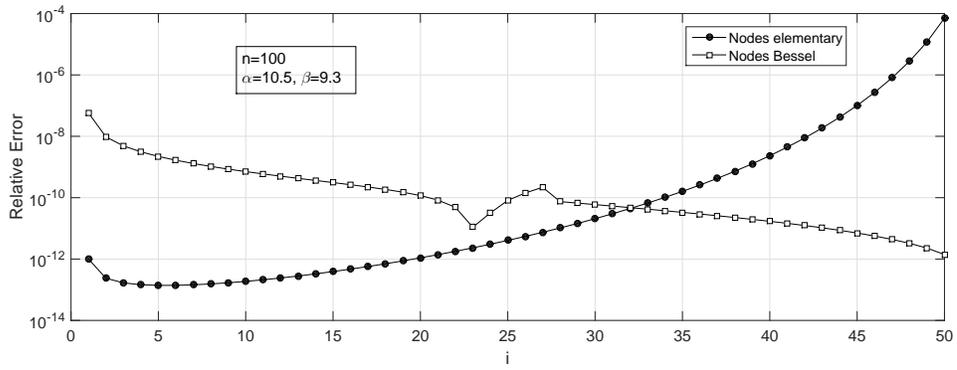}
\caption{
\label{pargra} Performance of the asymptotic expansions for computing the nodes
and scaled weights of $P^{(\alpha,\beta)}_n(x)$ for
 $\alpha=10.5$, $\beta=9.3$ and $n=1000$.
}
\end{figure}

Our algorithms are able to produce close to double precision accuracy for $n\ge 100$ and $-1<\alpha,\beta\le 5$. 
The accuracy worsens
as large parameters $\alpha$ and $\beta$ or smaller $n$ is considered. Larger values of the parameters $\alpha$ and $\beta$ 
are not very common in applications (in particular those related to 
barycentric interpolation and spectral methods which are discussed in the next section) 
because the Lebesgue constants becomes large, as discussed in  \cite{Wang:2014:EBW} (for more details 
\cite[sect. 14.4]{Szego:1975:OP}). As an illustration of a case with larger parameters, we shown in Fig. 
\ref{pargra}
the accuracy in the computation of the nodes for $n=100$, $\alpha=10.5$ and $\beta=9.3$; 
we observe that some accuracy is lost with respect
to the case $-1<\alpha,\beta\le 5$, but still the nodes can be computed with $10^{-10}$ relative accuracy.

\begin{figure}
\epsfxsize=15cm \epsfbox{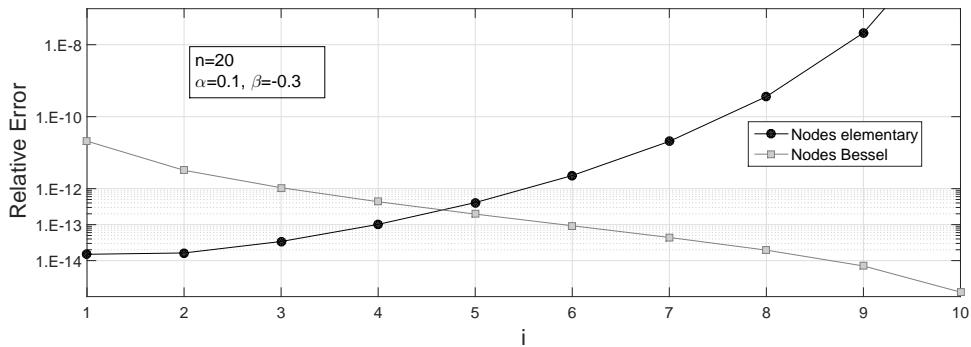}
\caption{
\label{n20} Performance of the asymptotic expansions for computing the nodes of $P^{(\alpha,\beta)}_n(x)$ for
 $\alpha=0.1$, $\beta=0.3$ and $n=20$.
}
\end{figure}

In the same way, our algorithms can also be used for smaller $n$ and they are quite accurate when $n\ge 20$ 
as long as the parameters are not large. For instance, as shown in Fig. \ref{n20} 
it is possible to obtain $10^{-12}$ relative accuracy for the nodes when $n=20$.  
This means that only one Newton step would be enough required to reach double precision accuracy in that cases,
thus improving the speed of iterative methods for small degrees.
   
\section{Extensions and further applications}
\label{paltonto}

In this paper, we have discussed the non-iterative computation of Gauss-Jacobi quadrature. As we next explain, the methods and techniques considered in this paper can also be applied in 
other contexts. 

\subsection{Computation of Jacobi polynomials}

First, we mention that the expansions we have described for Jacobi polynomials, together with the described
methods for obtaining the coefficients of the expansions, provide accurate methods for computing the Jacobi
polynomials for moderately large orders ($n\ge 100$) and $-1<\alpha,\beta\le 5$. 
We expect that the combination of these methods with the use of the three-term recurrence
relation for $n<100$, will lead to fast and accurate methods for computing Jacobi polynomials, similarly 
as happened for the computation of Laguerre polynomials in \cite{Gil:2017:ECL}.

\subsection{Jacobi-Gauss-Radau and Jacobi-Gauss-Lobatto quadratures}

In some applications (for example for solving boundary value problems by spectral methods 
\cite{Shen:2011:SM}), 
it can be interesting to include the boundaries of the integration interval as nodes,
maybe even with degree of multiplicity greater than one. The generalized Gauss-Lobatto rule in $[a,b]$ for
a weight function $w(x)$ supported in $[a,b]$ can be defined as the rule
$$
\displaystyle\int_{a}^{b} f(x) w(x) dx \approx
\displaystyle\sum_{j=0}^{j_{a}-1} v_0^{[j]} f^{(j)}(a)+\displaystyle\sum_{i=1}^{n} v_{i} f(x_i) +  
\displaystyle\sum_{j=0}^{j_{b}-1} v_{n+1}^{[j]} f^{(j)}(b),
$$ 
with the highest degree of accuracy, which is $2n+j_a+j_b-1$; 
here $f^{(j)}$ stands for the $j$-th derivative of $f$, and if a sum $\sum_{j=0}^{-1}$ appears we consider it
an empty sum. Generalized 
Jacobi-Gauss-Lobatto quadrature corresponds to $a=-1$, $b=+1$ and $w(x)=(1-x)^{\alpha}(1+x)^{\beta}$.
The most common cases in applications are when $j_a,j_b=1$ (Gauss-Lobatto) and $j_a=1,\,j_b=0$ or 
$j_a=0,\,j_b=1$ (Gauss-Radau).

For obvious reasons, the weights $v_{0}^{[j]}$ and $v_{n+1}^{[j]}$ are called boundary 
weights, while the values $x_{i}$ and $v_{i}$ are the internal nodes and weights. This
problem is closely related to the problem of Gaussian quadrature studied so far, and 
in particular the internal nodes and weights can be computed using exactly the same techniques,
because they are directly related to Gaussian quadrature. Indeed, it is easy to show that
the nodes $x_i$ are the Gaussian nodes in $[a,b]$ associated to the weight function 
$w(x)(x-a)^{j_a}(b-x)^{j_b}$, while the corresponding Gaussian weights $w_i$ for this modified weight function
are related with the internal weights of the generalized Gauss-Lobatto rule by 
$$v_i =\Frac{w_i}{(x_i-a)^{j_a}(b-x_i)^{j_b}}.$$ The computation of the boundary weights is a different problem,
but for the problem we are considering in this paper (the Jacobi case) there exist iterative 
formulas for computing them \cite{Petrova:2017:GGR}, and for the simple cases when $j_a,j_b\le 1$ (Gauss-Radau
or Gauss Lobatto), the boundary weights have simple and explicit expressions.

For instance, for the Jacobi-Gauss-Lobatto case we have 
$$
\displaystyle\int_{a}^{b} (1-x)^\alpha (1+x)^\beta f(x)  dx \approx
v_0 f(a)+\displaystyle\sum_{i=1}^{n} v_{i} f(x_i) +  v_{n+1} f(b),
$$ 
where $x_i$ are the zeros of $P_n^{(\alpha+1,\beta+1)} (x)$ and $v_i=w_i/(1-x_i^2)$, with $w_i$ the
weights for the $n$-degree Gauss-Jacobi quadrature with parameters $\alpha+1$ and $\beta+1$. Therefore, the
internal nodes and weights can be computed with the expressions given in this paper. For the boundary weights,
the following formula are available (see, for instance, \cite{Gautschi:2000:HOG}):
$$
v_0 = 2^{\alpha +\beta +1} \Frac{(\beta +1)\Gamma (\beta+1)^2 \Gamma (n+1) \Gamma (n+\alpha +2)}
{\Gamma (n+\beta+2) \Gamma (n+\alpha+\beta+3)},
$$
and $v_{n+1}$ has the same expression as $v_0$ but with $\alpha$ interchanged with $\beta$.

\subsection{Barycentric interpolation}

Let  $S=\{(x_i,f_i),\,i=1,\dots n\}$ be a set of points of ${\mathbb R}^2$, with $x_i\neq x_j$ if $i\neq j$. Then 
the polynomial of the smallest degree $Q(x)$ such that $Q(x_i)=f_i$, $i=1,\ldots n$ 
(the Lagrange interpolation polynomial for the set $S$) can be written as
\begin{equation}
\label{barye}
Q(x)=\Frac{\sum_{i=1}^n\Frac{u_i}{x-x_i}f_i}{\sum_{i=1}^n\Frac{u_i}{x-x_i}},\quad u_i=k
\left(\displaystyle\prod_{j\neq i}^n (x_i-x_j)\right)^{-1},
\end{equation}
with $k$ an arbitrary constant; the quantities $u_i$ are called barycentric weights.
This is called the barycentric formula of the Lagrange interpolation polynomial. Among the different 
formulas for computing Lagrange interpolants, this is one of the most efficient interpolation methods 
and it is a stable method \cite{Higham:2004:TNS}. 

For similar reasons that make the Gaussian quadrature optimal in the sense of the degree accuracy, 
it can be shown
that the interpolation at the orthogonal polynomial nodes have particularly good convergence rates 
(see for instance \cite{Xiang:2016:OIA}). A well-known example appearing in all textbooks 
is that of interpolation at the Chebyshev nodes, which smoothens the Runge phenomenon appearing when 
the nodes $\{x_i\}$ are equally spaced. We refer to \cite{Wang:2014:EBW} for additional examples of
interpolation using Jacobi nodes and weights (some of them requiring few hundreds of nodes for double
precision accuracy). In that paper the barycentric weights where computed with a Matlab version of the 
algorithm \cite{Glaser:2007:AFA}, and it is mentioned that the inclusion of the methods of 
\cite{Hale:2013:FAC} should improve the performance.

As discussed in \cite{Wang:2014:EBW}, and as it will be easy to understand for the reasons we will now explain,
the barycentric weights at the nodes of orthogonal polynomials can be easily  related to the Gauss 
weights for this set of nodes; the same is true for the Gauss-Radau and Gauss-Lobatto variants. To see this, 
first we observe that the barycentric weights of Eq. (\ref{barye}) can be written as (taking $k=1$):
$$
u_i (x)= 1/\phi'(x_i),\quad \phi (x)=\displaystyle\prod_{j=1}^{n}(x-x_j).
$$
Therefore, when the nodes are the roots of a polynomial (with only simple roots), the barycentric weights are inversely proportional to the derivative of the polynomial at the roots (with an arbitrary proportionality
constant). Then, the Gaussian weights are proportional to the square of the barycentric weights up to, 
possibly, an elementary function of the nodes (see Eq. (\ref{eq:numer1})). If we write the barycentric weights taking
the square root, the barycentric weights for interpolation at the zeros of $P_n^{(\alpha ,\beta)}(x)$ can
be written as
$$
u_i=(-1)^i \sin\theta_i \sqrt{w_i},
$$ 
where $\sin\theta_i=\sin(\arccos x_i)=\sqrt{1-x_i^2}$ and $w_i$ are the corresponding Gauss-Jacobi weights.
The alternating signs $(-1)^i$ are obvious because the derivative of the orthogonal polynomial takes 
different signs at consecutive zeros. Similar relations exist for the internal weights of 
Jacobi-Gauss-Radau ($x=-1$ or $x=+1$ are also nodes) and Jacobi-Gauss-Lobatto 
($x=-1$ or $x=+1$ are also nodes) interpolations, as can be understood from the previous discussion on the
associated quadrature rules.

Examples of applications employing barycentric interpolation at Jacobi points (or variants like 
Jacobi-Lobatto points) are the earlier mentioned references \cite{Wang:2014:EBW,Xiang:2016:OIA}. The case
of for Hermite-Fej\'er at Jacobi points was considered in \cite{Xiang:2015:TFI}.

\subsection{Other applications and some advantages with respect to previous methods}

Due to the many applications of Gauss quadrature and its variants (Gauss-Radau and Gauss-Lobatto) and
barycentric interpolation, it is not surprising that Gauss-Jacobi nodes and weights pop up in a many 
applications, for many of which a moderately large number of nodes may appear. 

For the particular case of Legendre quadratures, and as commented in 
\cite{Swarztrauber:2002:OCT}, atmospheric models use quadratures with a large number of nodes
(as an example of this see for instance \cite{Loeb:2002:TNS}, where 200 nodes Gauss-Legendre quadrature is employed). As mentioned in \cite{Hotta:2018:OP} (where quadratures with hundreds and even thousands of 
nodes are considered), 
one disadvantage of the Gaussian rules is that explicit formulas do not exists for
the nodes and weights and that iterative methods are usually needed; however, similarly as in 
\cite{Bogaert:2014:IFC}, we do 
provide explicit formulas which are iteration free, and that are fast enough to open the possibility of fast
computations with varying numbers of nodes.
 
But, of course, not only the Legendre case is important and also the Gauss-Jacobi rules appear 
in a vast number of applications, particularly for integrals with end-point singularities. Just to mention 
a few particular recent examples of the many which can be found in the literature see 
\cite{Yilmaz:2018:FRC} (mechanics), \cite{Heimlich:2016:PGI} (nuclear engineering) and 
\cite{Micheli:2013:TEI} (approximation theory).
On the other hand, the approximations in terms of Jacobi polynomials for different values of the
parameters (Legendre, Gegenbauer, Chebyshev or general Jacobi) have, as discussed in 
\cite{Boyd:2014:TRB}, different advantages depending on the criterion of merit, and different applications use 
some of these features. In most of these applications, the barycentric interpolation at the polynomial zeros
 is at the heart of the proposed methods. Without aiming at being exhaustive, apart from the applications 
mentioned in \cite{Bogaert:2014:IFC} for the particular case of Gauss-Legendre, we add the spectral 
Gegenbauer method of \cite{Elgindy:2017:HOS} (where hundreds of nodes and weights are used), the 
Jacobi pseudospectral methods applied to optimal control problems of reference \cite{Tang:2015:IFP} and fractional
 differential equations of \cite{Tang:2017:FPS}, for which the use of asymptotic methods allows 
for an application of the method for high accuracies (see Remark 14 of \cite{Tang:2015:IFP} and Remark 
4.5 of \cite{Tang:2017:FPS}) or the Jacobi spectral methods for the
Black-Scholes model of option pricing \cite{Nte:2014:ATB} which uses hundreds of nodes. 
Finally we recall that for the computation of the barycentric interpolation 
formulas, it is important to rely on methods which are efficient for large or moderately large degrees,
 both for Lagrange \cite{Wang:2014:EBW} as well as for Hermite-Fej\'er \cite{Xiang:2015:TFI} interpolation.

As we see, in a good number applications, and in particular for barycentric interpolation and
(pseudo-)spectral methods, a moderately large or even a vary large number of nodes and weights is needed. It 
would be completely inefficient (as illustrated in \cite{Hale:2013:FAC}) 
as well as unreliable to compute these nodes and weights by the Golub-Welsch
algorithm; instead, one should consider the more modern asymptotic methods, starting with  
\cite{Bogaert:2012:OCO,Hale:2013:FAC,Bogaert:2014:IFC,Townsend:2016:FCO}, continuing with our work 
on Gauss-Laguerre and Gauss-Hermite \cite{Gil:2018:GHL} and ending in the present paper, which extends the asymptotic
methods of Bogaert to the more general and important case of Jacobi, and has advantages both in terms 
of accuracy and speed with respect to \cite{Hale:2013:FAC}.

The present paper offers direct approximations for the Gauss-Jacobi nodes and weights in the same ranges as
\cite{Hale:2013:FAC} but in a faster non-iterative way, and with double precision relative accuracy both for the nodes
and the weights. We expect that these approximations will become part of 
forthcoming numerical software for computing Gauss-Jacobi nodes and weights and Jacobi barycentric weights.
It is also expected that the fast and direct computation of Gauss-Jacobi nodes will have an 
impact on the implementation of spectral methods for solving a huge variety of problems. The algorithms
are easily parallelizable, similarly as discussed in  \cite{Bogaert:2012:OCO,Bogaert:2014:IFC} for the
particular case of Gauss-Legendre,
and differently from iterative algorithms like \cite{Hale:2013:FAC} (which can not be parallelized so
 efficiently because each node requires several repeated iterations). 
Furthermore, in applications where just the nodes are needed the advantage over \cite{Hale:2013:FAC} will be 
even larger in terms of speed, because for that algorithm computing the nodes was almost as
expensive as computing the nodes and weights (as they require iterations), while our approximations for the 
nodes are much cheaper. When both the nodes and weights are computed, the speed-up will be in any case always 
larger than $2N$, with $N$ the number of iterations required by the iterative method; our numerical experiments
suggest that, in fact, the speed-up factor could be considerably larger. Also for smaller degrees $n$ ($n\le 20$) 
the approximations for the nodes are useful because they are more accurate than previous approximations and 
they can be used to speed-up iterative methods (just one Newton iteration would give more than double precision
accuracy).

\section*{Acknowledgements}
The authors thank the constructive comments of the reviewers. The authors also thank the associated editor for
suggesting that some applications of Gauss-Jacobi quadratures should be included in the paper.

\bibliographystyle{siam}
\bibliography{gauss}

\end{document}